\documentclass[10pt]{article}
\usepackage{amsmath}
\usepackage{amssymb}
\usepackage{amsthm}
\usepackage{mathrsfs}
\numberwithin{equation}{section}

\begin{document}
\title{Some estimates for commutators of fractional integrals associated to operators with Gaussian kernel bounds on weighted Morrey spaces}
\author{Hua Wang \footnote{E-mail address: wanghua@pku.edu.cn.
Supported by National Natural
Science Foundation of China under Grant \#10871173 and \#10931001.}\\
\footnotesize{Department of Mathematics, Zhejiang University, Hangzhou 310027, China}}
\date{}
\maketitle

\begin{abstract}
Let $L$ be the infinitesimal generator of an analytic semigroup on $L^2(\mathbb R^n)$ with Gaussian kernel bound, and let $L^{-\alpha/2}$ be the fractional integrals of $L$ for $0<\alpha<n$. In this paper, we will obtain some boundedness properties of commutators $\big[b,L^{-\alpha/2}\big]$ on the weighted Morrey spaces $L^{p,\kappa}(w)$ when the symbol $b$ belongs to $BMO(\mathbb R^n)$ or homogeneous Lipschitz space.
\end{abstract}

\section{Introduction and main results}

Suppose that $L$ is the infinitesimal generator of an analytic semigroup $\{e^{-tL}\}_{t>0}$ on $L^2(\mathbb R^n)$ with a kernel $p_t(x,y)$ satisfying a Gaussian upper bound; that is, there exist positive constants $C$ and $A$ such that for all $x$, $y\in\mathbb R^n$ and all $t>0$, we have
\begin{equation}
|p_t(x,y)|\le\frac{C}{t^{n/2}}e^{-A\frac{|x-y|^2}{t}}.
\end{equation}
Throughout this paper, we will always assume that the semigroup $\{e^{-tL}\}_{t>0}$ has a kernel which satisfies (1.1). This property is satisfied by a large class of differential operators, some examples can be seen in \cite{duong1}.

For any $0<\alpha<n$, the fractional integrals $L^{-\alpha/2}$ associated to the operator $L$ is defined by
\begin{equation}
L^{-\alpha/2}f(x)=\frac{1}{\Gamma(\alpha/2)}\int_0^\infty e^{-tL}(f)(x)t^{\alpha/2-1}\,dt.
\end{equation}
Note that if $L=-\Delta$ is the Laplacian on $\mathbb R^n$, then $L^{-\alpha/2}$ is the classical fractional integral operator $I_\alpha$, which is given by(see \cite{stein1})
\begin{equation*}
I_\alpha f(x)=\frac{\Gamma(\frac{n-\alpha}{2})}{2^{\alpha}\pi^{\frac n2}\Gamma(\frac{\alpha}2)}\int_{\mathbb R^n}\frac{f(y)}{|x-y|^{n-\alpha}}\,dy.
\end{equation*}
Let $b$ be a locally integrable function on $\mathbb R^n$. The commutator of $b$ and $L^{-\alpha/2}$ is defined as follows
\begin{equation}
\big[b,L^{-\alpha/2}\big](f)(x)=b(x)L^{-\alpha/2}(f)(x)-L^{-\alpha/2}(bf)(x).
\end{equation}

The first result on the theory of commutators was obtained by Coifman, Rochberg and Weiss in \cite{coifman}. Since then, many authors have been interested in studying this theory. When $0<\alpha<n$, $1<p<n/\alpha$ and $1/q=1/p-\alpha/n$, Chanillo \cite{chanillo} proved that the commutator $[b,I_\alpha]$ is bounded from $L^p(\mathbb R^n)$ to $L^q(\mathbb R^n)$ whenever $b\in BMO(\mathbb R^n)$. Paluszy\'{n}ski \cite{paluszynski} showed that $b\in {\dot\Lambda}_\beta(\mathbb R^n)$ (homogeneous Lipschitz space) if and only if $[b,I_\alpha]$ is bounded from $L^p(\mathbb R^n)$ to $L^s(\mathbb R^n)$, where $0<\beta<1$, $1<p<n/{(\alpha+\beta)}$ and $1/s=1/p-{(\alpha+\beta)}/n$. For the weighted case, Segovia and Torrea \cite{segovia} proved that when $b\in BMO(\mathbb R^n)$ and $w\in A_{p,q}$(Muckenhoupt weight class), then $[b,I_\alpha]$ is bounded from $L^p(w^p)$ to $L^q(w^q)$.

In 2004, by using a new sharp maximal function introduced in \cite{martell}, Duong and Yan \cite{duong1} extended the result of \cite{chanillo} from $(-\Delta)$ to the more general operator $L$ defined above. More precisely, they showed that

\newtheorem*{thma}{Theorem A}
\begin{thma}
Let $0<\alpha<n$, $1<p<n/\alpha$ and $1/q=1/p-\alpha/n$. If $b\in BMO(\mathbb R^n)$, then the commutator $\big[b,L^{-\alpha/2}\big]$ is bounded from $L^p(\mathbb R^n)$ to $L^q(\mathbb R^n)$.
\end{thma}

In 2008, Auscher and Martell \cite{auscher} considered the weighted case and obtained the following result (see also \cite{cruz}).

\newtheorem*{thmb}{Theorem B}
\begin{thmb}
Let $0<\alpha<n$, $1<p<n/\alpha$, $1/q=1/p-\alpha/n$ and $w\in A_{p,q}$. If $b\in BMO(\mathbb R^n)$, then the commutator $\big[b,L^{-\alpha/2}\big]$ is bounded from $L^p(w^p)$ to $L^q(w^q)$.
\end{thmb}

On the other hand, in 2009, Komori and Shirai \cite{komori} first introduced the weighted Morrey spaces $L^{p,\kappa}(w)$ which could be viewed as an extension of weighted Lebesgue spaces, and investigated the boundedness of the Hardy-Littlewood maximal operator, singular integral operator and fractional integral operator on these weighted spaces. Moreover, they also proved the following theorem.

\newtheorem*{thmc}{Theorem C}
\begin{thmc}
Let $0<\alpha<n$, $1<p<n/\alpha$, $1/q=1/p-\alpha/n$, $0<\kappa<p/q$ and $w\in A_{p,q}$. If $b\in BMO(\mathbb R^n)$, then the commutator $[b,I_\alpha]$ is bounded from $L^{p,\kappa}(w^p,w^q)$ to $L^{q,{\kappa q}/p}(w^q)$.
\end{thmc}

The purpose of this paper is to study the boundedness of $\big[b,L^{-\alpha/2}\big]$ on the weighted Morrey spaces $L^{p,\kappa}(w)$ when $b\in BMO(\mathbb R^n)$ or $b\in{\dot\Lambda}_\beta(\mathbb R^n)$. Our main results in the paper are formulated as follows.

\newtheorem{theorem}{Theorem}[section]

\begin{theorem}
Let $0<\alpha<n$, $1<p<n/\alpha$, $1/q=1/p-\alpha/n$, $0<\kappa<p/q$ and $w\in A_{p,q}$. Suppose that $b\in BMO(\mathbb R^n)$, then the commutator $\big[b,L^{-\alpha/2}\big]$ is bounded from $L^{p,\kappa}(w^p,w^q)$ to $L^{q,{\kappa q}/p}(w^q)$.
\end{theorem}

\begin{theorem}
Let $0<\beta<1$, $0<\alpha+\beta<n$, $1<p<n/{(\alpha+\beta)}$, $1/s=1/p-{(\alpha+\beta)}/n$ and $w\in A_{p,s}$. Suppose that $b\in{\dot\Lambda}_\beta(\mathbb R^n)$ and $0<\kappa<p/s$, then the commutator $\big[b,L^{-\alpha/2}\big]$ is bounded from $L^{p,\kappa}(w^p,w^s)$ to $L^{s,{\kappa s}/p}(w^s)$.
\end{theorem}

\begin{theorem}
Let $0<\beta<1$, $0<\alpha+\beta<n$, $1/s=1-{(\alpha+\beta)}/n$ and $w\in A_{1,s}$. Suppose that $b\in{\dot\Lambda}_\beta(\mathbb R^n)$ and $0<\kappa<1/s$, then the commutator $\big[b,L^{-\alpha/2}\big]$ is bounded from $L^{1,\kappa}(w,w^s)$ to $WL^{s,{\kappa s}}(w^s)$.
\end{theorem}

\begin{theorem}
Let $0<\beta<1$, $0<\alpha+\beta<n$, $1<p<n/{(\alpha+\beta)}$, $1/s=1/p-{(\alpha+\beta)}/n$ and $w\in A_{p,s}$. Suppose that $b\in{\dot\Lambda}_\beta(\mathbb R^n)$ and $\kappa=p/s$, then the commutator $\big[b,L^{-\alpha/2}\big]$ is bounded from $L^{p,\kappa}(w^p,w^s)$ to $BMO_L$.
\end{theorem}

\section{Notations and definitions}

First let us recall some standard definitions and notations. The classical $A_p$ weight theory was first introduced by Muckenhoupt in the study of weighted $L^p$ boundedness of Hardy-Littlewood maximal functions in \cite{muckenhoupt1}. Let $w$ be a nonnegative, locally integrable function defined on $\mathbb R^n$, $B=B(x_0,r_B)$ denotes the ball with the center $x_0$ and radius $r_B$. Given a ball $B$ and $\lambda>0$, $\lambda B$ denotes the ball with the same center as $B$ whose radius is $\lambda$ times that of $B$, we also denote the Lebesgue measure of $B$ by $|B|$ and the weighted measure of $B$ by $w(B)$, where $w(B)=\int_B w(x)\,dx$.
We say that $w\in A_p$, $1<p<\infty$, if
\begin{equation*}
\left(\frac1{|B|}\int_B w(x)\,dx\right)\left(\frac1{|B|}\int_B w(x)^{-1/{(p-1)}}\,dx\right)^{p-1}\le C \quad\mbox{for every ball}\; B\subseteq \mathbb
R^n,
\end{equation*}
where $C$ is a positive constant which is independent of $B$.

For the case $p=1$, $w\in A_1$, if
\begin{equation*}
\frac1{|B|}\int_B w(x)\,dx\le C\cdot\underset{x\in B}{\mbox{ess\,inf}}\,w(x)\quad\mbox{for every ball}\;B\subseteq\mathbb R^n.
\end{equation*}

For the case $p=\infty$, $w\in A_\infty$ if it satisfies the $A_p$ condition for some $1<p<\infty$.

We also need another weight class $A_{p,q}$ introduced by Muckenhoupt and Wheeden in \cite{muckenhoupt2}. A weight function $w$ belongs to $A_{p,q}$ for $1<p<q<\infty$ if there exists a constant $C>0$ such that
\begin{equation*}
\left(\frac{1}{|B|}\int_B w(x)^q\,dx\right)^{1/q}\left(\frac{1}{|B|}\int_B w(x)^{-p'}\,dx\right)^{1/{p'}}\le C \quad\mbox{for every ball}\; B\subseteq \mathbb R^n,
\end{equation*}
where $p'$ denotes the conjugate exponent of $p>1$; that is, $1/p+1/{p'}=1$.

When $p=1$, $w$ is in the class $A_{1,q}$ with $1<q<\infty$ if there exists a constant $C>0$ such that
\begin{equation*}
\left(\frac{1}{|B|}\int_B w(x)^q\,dx\right)^{1/q}\bigg(\underset{x\in B}{\mbox{ess\,sup}}\,\frac{1}{w(x)}\bigg)\le C \quad\mbox{for every ball}\; B\subseteq \mathbb R^n.
\end{equation*}

A weight function $w$ is said to belong to the reverse H\"{o}lder class $RH_r$ if there exist two constants $r>1$ and $C>0$ such that the following reverse H\"{o}lder inequality holds
\begin{equation*}
\left(\frac{1}{|B|}\int_B w(x)^r\,dx\right)^{1/r}\le C\left(\frac{1}{|B|}\int_B w(x)\,dx\right)\quad\mbox{for every ball}\; B\subseteq \mathbb R^n.
\end{equation*}

We give the following results that we will use in the sequel.

\newtheorem{lemma}[theorem]{Lemma}
\begin{lemma}[\cite{garcia}]
Let $w\in A_p$ with $p\ge1$. Then, for any ball $B$, there exists an absolute constant $C>0$ such that
\begin{equation*}
w(2B)\le C\,w(B).
\end{equation*}
In general, for any $\lambda>1$, we have
\begin{equation*}
w(\lambda B)\le C\cdot\lambda^{np}w(B),
\end{equation*}
where $C$ does not depend on $B$ nor on $\lambda$.
\end{lemma}

\begin{lemma}[\cite{gundy}]
Let $w\in RH_r$ with $r>1$. Then there exists a constant $C>0$ such that
\begin{equation*}
\frac{w(E)}{w(B)}\le C\left(\frac{|E|}{|B|}\right)^{(r-1)/r}
\end{equation*}
for any measurable subset $E$ of a ball $B$.
\end{lemma}

Next we shall give the definitions of some function spaces. A locally integrable function $b$ is said to be in $BMO(\mathbb R^n)$ if
\begin{equation*}
\|b\|_*=\sup_{B}\frac{1}{|B|}\int_B|b(x)-b_B|\,dx<\infty,
\end{equation*}
where $b_B=\frac{1}{|B|}\int_B b(y)\,dy$ and the supremum is taken over all balls $B$ in $\mathbb R^n$.

\begin{theorem}[\cite{duoand,john}]
Assume that $b\in BMO(\mathbb R^n)$. Then for any $1\le p<\infty$, we have
\begin{equation*}
\sup_B\bigg(\frac{1}{|B|}\int_B\big|b(x)-b_B\big|^p\,dx\bigg)^{1/p}\le C\|b\|_*.
\end{equation*}
\end{theorem}

Let $0<\beta\le1$. The homogeneous Lipschitz space ${\dot\Lambda}_\beta(\mathbb R^n)$ is defined by
\begin{equation*}
{\dot\Lambda}_\beta(\mathbb R^n)=\Big\{b:\|b\|_{{\dot\Lambda}_\beta}=\sup_{x,h\in\mathbb R^n,h\neq0}\frac{|b(x+h)-b(x)|}{|h|^{\beta}}<\infty\Big\}.
\end{equation*}

Given a Muckenhoupt's weight function $w$ on $\mathbb R^n$, for $1\le p<\infty$, we denote by $L^p(w)$ the space of all functions $f$ satisfying
\begin{equation*}
\|f\|_{L^p(w)}=\bigg(\int_{\mathbb R^n}|f(x)|^pw(x)\,dx\bigg)^{1/p}<\infty.
\end{equation*}
When $p=\infty$, $L^\infty(w)$ will be taken to mean $L^\infty(\mathbb R^n)$, and
\begin{equation*}
\|f\|_{L^\infty(w)}=\|f\|_{L^\infty}=\underset{x\in\mathbb R^n}{\mbox{ess\,sup}}\,|f(x)|.
\end{equation*}

In \cite{komori}, Komori and Shirai first defined the weighted Morrey spaces and obtained some known results relevant to this paper. For the boundedness of some other operators on these spaces, we refer the readers to \cite{wang1,wang2}.

\newtheorem{defn}[theorem]{Definition}
\begin{defn}[\cite{komori}]
Let $1\le p<\infty$, $0<\kappa<1$ and $w$ be a weight function. Then the weighted Morrey space is defined by
\begin{equation*}
L^{p,\kappa}(w)=\big\{f\in L^p_{loc}(w):\|f\|_{L^{p,\kappa}(w)}<\infty\big\},
\end{equation*}
where
\begin{equation*}
\|f\|_{L^{p,\kappa}(w)}=\sup_B\left(\frac{1}{w(B)^\kappa}\int_B|f(x)|^pw(x)\,dx\right)^{1/p}
\end{equation*}
and the supremum is taken over all balls $B$ in $\mathbb R^n$.
\end{defn}

We also denote by $WL^{p,\kappa}(w)$ the weighted weak Morrey space of all locally integrable functions satisfying
\begin{equation*}
\|f\|_{WL^{p,\kappa}(w)}=\sup_B\sup_{t>0}\frac{1}{w(B)^{\kappa/p}}t\cdot w\big(\big\{x\in B:|f(x)|>t\big\}\big)^{1/p}<\infty.
\end{equation*}

In order to deal with the fractional order case, we need to consider the weighted Morrey space with two weights.

\begin{defn}[\cite{komori}]
Let $1\le p<\infty$ and $0<\kappa<1$. Then for two weights $u$ and $v$, the weighted Morrey space is defined by
\begin{equation*}
L^{p,\kappa}(u,v)=\big\{f\in L^p_{loc}(u):\|f\|_{L^{p,\kappa}(u,v)}<\infty\big\},
\end{equation*}
where
\begin{equation*}
\|f\|_{L^{p,\kappa}(u,v)}=\sup_{B}\left(\frac{1}{v(B)^{\kappa}}\int_B|f(x)|^pu(x)\,dx\right)^{1/p}.
\end{equation*}
\end{defn}

\newtheorem*{theoremd}{Theorem D}
\begin{theoremd}[\cite{komori}]
Let $0<\alpha<n$, $1<p<n/\alpha$, $1/q=1/p-\alpha/n$ and $w\in A_{p,q}$. Then the operator $I_\alpha$ is bounded from $L^{p,\kappa}(w^p,w^q)$ to $L^{q,{\kappa q}/p}(w^q)$.
\end{theoremd}

\newtheorem*{theoreme}{Theorem E}
\begin{theoreme}[\cite{komori}]
Let $p=1$, $0<\alpha<n$, $1/q=1-\alpha/n$ and $w\in A_{1,q}$. Then the operator $I_\alpha$ is bounded from $L^{1,\kappa}(w,w^q)$ to $WL^{q,{\kappa q}}(w^q)$.
\end{theoreme}

We are going to conclude this section by defining the function spaces $BMO_L$. Duong and Yan \cite{duong3} introduced and developed a new function space $BMO_L$ associated with an operator $L$. Assume that the kernel $p_t(x,y)$ of $\{e^{-tL}\}_{t>0}$ satisfies an upper bound
\begin{equation*}
|p_t(x,y)|\le t^{-n/2}g\Big(\frac{|x-y|}{t^{1/2}}\Big),
\end{equation*}
for all $x,y\in\mathbb R^n$ and all $t>0$. Here $g$ is a positive, bounded, decreasing function satisfying
\begin{equation}
\lim_{r\to\infty}r^{n+\epsilon}g(r)=0, \quad \mbox{for some}\; \epsilon>0.
\end{equation}
Let $\epsilon$ be the constant in (2.1) and $0<\beta<\epsilon$. A function $f\in L^p_{loc}(\mathbb R^n)$ is said to be a function of type $(p,\beta)$ if $f$ satisfies
\begin{equation}
\bigg(\int_{\mathbb R^n}\frac{|f(x)|^p}{(1+|x|)^{n+\beta}}\,dx\bigg)^{1/p}\le c<\infty.
\end{equation}
We denote by $\mathcal M_{(p,\beta)}$ the collection of all functions of type $(p,\beta)$. If $f\in\mathcal M_{(p,\beta)}$, then the norm of $f$ in $\mathcal M_{(p,\beta)}$ is defined by
\begin{equation*}
\|f\|_{\mathcal M_{(p,\beta)}}=\inf\big\{c\ge0:\mbox{(2.2) holds}\big\}.
\end{equation*}
It is easy to see that $\mathcal M_{(p,\beta)}$ is a Banach space under the norm $\|f\|_{\mathcal M_{(p,\beta)}}<\infty$. We set
\begin{equation*}
\mathcal M_p=\bigcup_{\beta:0<\beta<\epsilon}\mathcal M_{(p,\beta)}.
\end{equation*}
For any $f\in L^p(\mathbb R^n)$, $1\le p<\infty$, Martell \cite{martell} defined a kind of sharp maximal function $M^{\#}_L f$ associated with the semigroup $\{e^{-tL}\}_{t>0}$ by the expression
$$M^{\#}_L f(x)=\sup_{x\in B}\frac{1}{|B|}\int_B\left|f(y)-e^{-t_B L}f(y)\right|\,dy,$$
where $t_B=r_B^2$ and $r_B$ is the radius of the ball $B$. Let $f\in\mathcal M_p$ with $1<p<\infty$, then we say that $f\in BMO_L$ if the sharp maximal function $M^{\#}_L f\in L^\infty(\mathbb R^n)$, and we define $\|f\|_{BMO_L}=\big\|M^{\#}_L f\big\|_{L^\infty}$. For further details about the properties and applications of $BMO_L$ spaces, we refer the readers to \cite{deng,duong3,duong2}.

Throughout this article, we will use $C$ to denote a positive constant, which is independent of the main parameters and not necessarily the same at each occurrence. By $A\sim B$, we mean that there exists a constant $C>1$ such that $\frac1C\le\frac AB\le C$. Moreover, we denote the conjugate exponent of $q>1$ by $q'=q/{(q-1)}$.

\section{Proofs of Theorems 1.1, 1.2 and 1.3}

\begin{proof}
Fix a ball $B=B(x_0,r_B)\subseteq\mathbb R^n$ and decompose $f=f_1+f_2$, where $f_1=f\chi_{_{2B}}$, $\chi_{_{2B}}$ denotes the characteristic function of $2B$. Since $\big[b,L^{-\alpha/2}\big]$ is a linear operator, then we have
\begin{equation*}
\begin{split}
&\frac{1}{w^q(B)^{\kappa/p}}\bigg(\int_B\big|\big[b,L^{-\alpha/2}\big]f(x)\big|^qw^q(x)\,dx\bigg)^{1/q}\\
\le&\,\frac{1}{w^q(B)^{\kappa/p}}\bigg(\int_B\big|\big[b,L^{-\alpha/2}\big]f_1(x)\big|^qw^q(x)\,dx\bigg)^{1/q}
+\frac{1}{w^q(B)^{\kappa/p}}\bigg(\int_B\big|\big[b,L^{-\alpha/2}\big]f_2(x)\big|^qw^q(x)\,dx\bigg)^{1/q}\\
=&\,I_1+I_2.
\end{split}
\end{equation*}
For the term $I_1$, since $w\in A_{p,q}$, then we get $w^q\in A_{1+q/{p'}}$(see \cite{muckenhoupt2}). Hence, it follows from Theorem B and Lemma 2.1 that
\begin{align}
I_1&\le\frac{1}{w^q(B)^{\kappa/p}}\big\|\big[b,L^{-\alpha/2}\big]f_1\big\|_{L^q(w^q)}\notag\\
&\le C\|b\|_*\cdot\frac{1}{w^q(B)^{\kappa/p}}\|f_1\|_{L^p(w^p)}\notag\\
&=C\|b\|_*\cdot\frac{1}{w^q(B)^{\kappa/p}}\bigg(\int_{2B}|f(x)|^pw^p(x)\,dx\bigg)^{1/p}\notag\\
&\le C\|b\|_*\|f\|_{L^{p,\kappa}(w^p,w^q)}\cdot\frac{w^q(2B)^{\kappa/p}}{w^q(B)^{\kappa/p}}\notag\\
&\le C\|b\|_*\|f\|_{L^{p,\kappa}(w^p,w^q)}.
\end{align}
We now turn to deal with the term $I_2$. Denote the kernel of $L^{-\alpha/2}$ by $K_\alpha(x,y)$, then for any $x\in B$, we write
\begin{equation*}
\begin{split}
\big|\big[b,L^{-\alpha/2}\big]f_2(x)\big|=&\,\bigg|\int_{(2B)^c}\big[b(x)-b(y)\big]K_\alpha(x,y)f(y)\,dy\bigg|\\
\le&\,\big|b(x)-b_B\big|\cdot\int_{(2B)^c}|K_\alpha(x,y)||f(y)|\,dy\\
&+\int_{(2B)^c}|b(y)-b_B||K_\alpha(x,y)||f(y)|\,dy\\
=&\,\mbox{\upshape I+II}.
\end{split}
\end{equation*}
Since the kernel of $e^{-tL}$ is $p_t(x,y)$, then it follows immediately from (1.2) that (see \cite{mo})
\begin{equation}
K_\alpha(x,y)=\frac{1}{\Gamma(\alpha/2)}\int_0^\infty p_t(x,y)t^{\alpha/2-1}\,dt.
\end{equation}
Thus, by using the Gaussian upper bound (1.1) and the expression (3.2), we can deduce (see \cite{duong1} and \cite{mo})
\begin{align}
|K_\alpha(x,y)|&\le \frac{1}{\Gamma(\alpha/2)}\int_0^\infty |p_t(x,y)|t^{\alpha/2-1}\,dt\notag\\
&\le C\int_0^\infty e^{-A\frac{|x-y|^2}{t}}t^{\alpha/2-n/2-1}\,dt\notag\\
&\le C\cdot\frac{1}{|x-y|^{n-\alpha}}.
\end{align}
So we have
\begin{equation*}
\mbox{\upshape I}\le \big|b(x)-b_B\big|\cdot\sum_{k=1}^\infty\frac{1}{|2^{k+1}B|^{1-\alpha/n}}\int_{2^{k+1}B}|f(y)|\,dy.
\end{equation*}
By using H\"{o}lder's inequality and the fact that $w\in A_{p,q}$, we can get
\begin{align}
\int_{2^{k+1}B}|f(y)|\,dy&\le\bigg(\int_{2^{k+1}B}|f(y)|^pw^p(y)\,dy\bigg)^{1/p}
\bigg(\int_{2^{k+1}B}w^{-p'}(y)\,dy\bigg)^{1/{p'}}\notag\\
&\le C\|f\|_{L^{p,\kappa}(w^p,w^q)}\big|2^{k+1}B\big|^{1/q+1/{p'}}\cdot\frac{1}{w^q(2^{k+1}B)^{1/q-\kappa/p}}\notag\\
&= C\|f\|_{L^{p,\kappa}(w^p,w^q)}\big|2^{k+1}B\big|^{1-\alpha/n}\cdot\frac{1}{w^q(2^{k+1}B)^{1/q-\kappa/p}}.
\end{align}
Hence
\begin{equation*}
\begin{split}
&\frac{1}{w^q(B)^{\kappa/p}}\bigg(\int_B \mbox{\upshape I}^q\,w^q(x)\,dx\bigg)^{1/q}\\
\le&\, C\|f\|_{L^{p,\kappa}(w^p,w^q)}\frac{1}{w^q(B)^{\kappa/p}}
\sum_{k=1}^\infty\frac{1}{w^q(2^{k+1}B)^{1/q-\kappa/p}}\cdot\bigg(\int_B|b(x)-b_B|^qw^q(x)\,dx\bigg)^{1/q}\\
=&\, C\|f\|_{L^{p,\kappa}(w^p,w^q)}\sum_{k=1}^\infty\frac{w^q(B)^{1/q-\kappa/p}}{w^q(2^{k+1}B)^{1/q-\kappa/p}}
\cdot\bigg(\frac{1}{w^q(B)}\int_B|b(x)-b_B|^qw^q(x)\,dx\bigg)^{1/q}.
\end{split}
\end{equation*}
We now claim that for any $1<q<\infty$ and $v\in A_\infty$, the following inequality holds
\begin{equation}
\bigg(\frac{1}{v(B)}\int_B|b(x)-b_B|^qv(x)\,dx\bigg)^{1/q}\le C\|b\|_*.
\end{equation}
In fact, since $v\in A_\infty$, then we know that there exists $r>1$ such that $v\in RH_r$. Thus, by H\"older's inequality and Theorem 2.3, we obtain
\begin{equation*}
\begin{split}
\bigg(\frac{1}{v(B)}\int_B|b(x)-b_B|^qv(x)\,dx\bigg)^{1/q}&\le\frac{1}{v(B)^{1/q}}
\bigg(\int_B|b(x)-b_B|^{qr'}\,dx\bigg)^{1/{(qr')}}\bigg(\int_B v(x)^r\,dx\bigg)^{1/{(qr)}}\\
&\le C\bigg(\frac{1}{|B|}\int_B|b(x)-b_B|^{qr'}\,dx\bigg)^{1/{(qr')}}\\
&\le C\|b\|_*,
\end{split}
\end{equation*}
which is our desired result. Note that $w^q\in A_{1+q/{p'}}\subset A_\infty$. In addition, we have $w^q\in RH_r$ with $r>1$. Thus, by Lemma 2.2, we get
\begin{equation}
\frac{w^q(B)}{w^q(2^{k+1}B)}\le C\left(\frac{|B|}{|2^{k+1}B|}\right)^{{(r-1)}/r}.
\end{equation}
Consequently
\begin{align}
\frac{1}{w^q(B)^{\kappa/p}}\bigg(\int_B \mbox{\upshape I}^q\,w^q(x)\,dx\bigg)^{1/q}&\le C\|b\|_*\|f\|_{L^{p,\kappa}(w^p,w^q)}\sum_{k=1}^\infty\left(\frac{1}{2^{kn}}\right)^{(1-1/r)(1/q-\kappa/p)}\notag\\
&\le C\|b\|_*\|f\|_{L^{p,\kappa}(w^p,w^q)},
\end{align}
where the last series is convergent since $r>1$ and $0<\kappa<p/q$. On the other hand
\begin{equation*}
\begin{split}
\mbox{\upshape II}\le&\,\sum_{k=1}^\infty\int_{2^{k+1}B\backslash 2^kB}|b(y)-b_B||K_\alpha(x,y)||f(y)|\,dy\\
\le&\,\sum_{k=1}^\infty\int_{2^{k+1}B\backslash 2^kB}\big|b(y)-b_{2^{k+1}B}\big||K_\alpha(x,y)||f(y)|\,dy\\
&+\sum_{k=1}^\infty\int_{2^{k+1}B\backslash 2^kB}\big|b_{2^{k+1}B}-b_B\big||K_\alpha(x,y)||f(y)|\,dy\\
=&\,\mbox{\upshape III+IV}.
\end{split}
\end{equation*}
To estimate III and IV, we observe that when $x\in B$, $y\in (2B)^c$, then $|y-x|\sim|y-x_0|$. Hence, it follows directly from the kernel estimate (3.3) that
\begin{equation*}
\mbox{\upshape III}\le C \sum_{k=1}^\infty\frac{1}{|2^{k+1}B|^{1-\alpha/n}}
\int_{2^{k+1}B}\big|b(y)-b_{2^{k+1}B}\big||f(y)|\,dy.
\end{equation*}
An application of H\"older's inequality yields
\begin{equation*}
\begin{split}
&\int_{2^{k+1}B}\big|b(y)-b_{2^{k+1}B}\big||f(y)|\,dy\\
\le&\,
\bigg(\int_{2^{k+1}B}\big|b(y)-b_{2^{k+1}B}\big|^{p'}w^{-p'}(y)\,dy\bigg)^{1/{p'}}
\bigg(\int_{2^{k+1}B}|f(y)|^pw^p(y)\,dy\bigg)^{1/p}\\
\le&\,\|f\|_{L^{p,\kappa}(w^p,w^q)}\cdot w^q\big(2^{k+1}B\big)^{\kappa/p}
\bigg(\int_{2^{k+1}B}\big|b(y)-b_{2^{k+1}B}\big|^{p'}w^{-p'}(y)\,dy\bigg)^{1/{p'}}.
\end{split}
\end{equation*}
We set $v(y)=w^{-p'}(y)$, then we have $v\in A_{1+{p'}/q}\subset A_\infty$ because $w\in A_{p,q}$ (see \cite{muckenhoupt2}). By the previous estimate (3.5) and the fact that $w\in A_{p,q}$, we obtain
\begin{align}
\bigg(\int_{2^{k+1}B}\big|b(y)-b_{2^{k+1}B}\big|^{p'}v(y)\,dy\bigg)^{1/{p'}}&\le C\|b\|_*v\big(2^{k+1}B\big)^{1/{p'}}\notag\\
&\le C\|b\|_*\cdot\frac{|2^{k+1}B|^{1/q+1/{p'}}}{w^q(2^{k+1}B)^{1/q}}.
\end{align}
Note that $1/q+1/{p'}=1-\alpha/n$. Hence, by (3.6) and (3.8), we have
\begin{align}
\frac{1}{w^q(B)^{\kappa/p}}\bigg(\int_B \mbox{\upshape III}^q\,w^q(x)\,dx\bigg)^{1/q}
&\le C\|b\|_*\|f\|_{L^{p,\kappa}(w^p,w^q)}
\sum_{k=1}^\infty\frac{w^q(B)^{1/q-\kappa/p}}{w^q(2^{k+1}B)^{1/q-\kappa/p}}\notag\\
&\le C\|b\|_*\|f\|_{L^{p,\kappa}(w^p,w^q)}.
\end{align}
Since $b\in BMO(\mathbb R^n)$, then a simple calculation gives that
\begin{equation*}
\big|b_{2^{k+1}B}-b_B\big|\le C\cdot k\|b\|_*.
\end{equation*}
Thus, by the estimates (3.3) and (3.4), we get
\begin{equation*}
\begin{split}
\mbox{\upshape IV}&\le C\|b\|_*\sum_{k=1}^\infty k\cdot\frac{1}{|2^{k+1}B|^{1-\alpha/n}}
\int_{2^{k+1}B}|f(y)|\,dy \\
&\le C\|b\|_*\|f\|_{L^{p,\kappa}(w^p,w^q)}\sum_{k=1}^\infty k\cdot\frac{1}{w^q(2^{k+1}B)^{1/q-\kappa/p}}.
\end{split}
\end{equation*}
Therefore
\begin{align}
\frac{1}{w^q(B)^{\kappa/p}}\bigg(\int_B \mbox{\upshape IV}^q\,w^q(x)\,dx\bigg)^{1/q}
&\le C\|b\|_*\|f\|_{L^{p,\kappa}(w^p,w^q)}
\sum_{k=1}^\infty k\cdot\frac{w^q(B)^{1/q-\kappa/p}}{w^q(2^{k+1}B)^{1/q-\kappa/p}}\notag\\
&\le C\|b\|_*\|f\|_{L^{p,\kappa}(w^p,w^q)}\sum_{k=1}^\infty \frac{k}{2^{kn\delta}}\notag\\
&\le C\|b\|_*\|f\|_{L^{p,\kappa}(w^p,w^q)},
\end{align}
where $w^q\in RH_r$ and $\delta=(1-1/r)(1/q-\kappa/p)$. Summarizing the estimates (3.9)
and (3.10) derived above, we thus obtain
\begin{equation}
\frac{1}{w^q(B)^{\kappa/p}}\bigg(\int_B \mbox{\upshape II}^q\,w^q(x)\,dx\bigg)^{1/q}\le C\|b\|_*\|f\|_{L^{p,\kappa}(w^p,w^q)}.
\end{equation}
Combining the inequalities (3.1) and (3.7) with the above inequality (3.11) and taking the supremum over all balls $B\subseteq\mathbb R^n$, we complete the proof of Theorem 1.1.
\end{proof}

Obviously, by (3.3), we have the following pointwise inequality
\begin{equation*}
\big|L^{-\alpha/2}(f)(x)\big|\le C\cdot I_\alpha(|f|)(x) \quad\mbox{for all}\;x\in{\mathbb R^n}.
\end{equation*}
Furthermore, by the definition of $b\in{\dot\Lambda}_\beta(\mathbb R^n)$ and (3.3), we deduce
\begin{align}
\big|\big[b,L^{-\alpha/2}\big](f)(x)\big|&\le\int_{\mathbb R^n}|b(x)-b(y)||K_\alpha(x,y)||f(y)|\,dy\notag\\
&\le C\|b\|_{{\dot\Lambda}_\beta}\int_{\mathbb R^n}\frac{|f(y)|}{|x-y|^{n-\alpha-\beta}}\,dy\notag\\
&\le C\|b\|_{{\dot\Lambda}_\beta}I_{\alpha+\beta}(|f|)(x).
\end{align}
Hence, Theorems 1.2 and 1.3 follows immediately from Theorems D and E.

\section{Proof of Theorem 1.4}

\begin{proof}
For any given $x\in \mathbb R^n$, fix a ball $B=B(x_0,r_B)$ which contains $x$. We decompose $f=f_1+f_2$, where $f_1=f\chi_{2B}$, and set $t_B=r^2_B$. Then we write
\begin{equation*}
\begin{split}
&\frac{1}{|B|}\int_B\left|\big[b,L^{-\alpha/2}\big]f(y)-e^{-t_B L}\big[b,L^{-\alpha/2}\big]f(y)\right|\,dy\\
\le&\,\frac{1}{|B|}\int_B\left|\big[b,L^{-\alpha/2}\big]f_1(y)\right|\,dy+\frac{1}{|B|}\int_B
\left|e^{-t_BL}\big[b,L^{-\alpha/2}\big]f_1(y)\right|\,dy\\
&+\frac{1}{|B|}\int_B\left|\big[b,L^{-\alpha/2}\big]f_2(y)-e^{-t_BL}\big[b,L^{-\alpha/2}\big]f_2(y)\right|\,dy\\
=&\,J_1+J_2+J_3.
\end{split}
\end{equation*}
We are now going to estimate each term respectively. For the first term $J_1$, since $w\in A_{p,s}$, then the operator $I_{\alpha+\beta}$ is bounded from $L^p(w^p)$ into $L^s(w^s)$(see \cite{muckenhoupt2}). We also know that $w^s\in A_{1+s/{p'}}\subset A_s$. Applying H\"{o}lder's inequality, the inequality (3.12), Lemma 2.1 and the fact that $w^s\in A_s$, we obtain
\begin{equation*}
\begin{split}
J_1&\le\frac{1}{|B|}\bigg(\int_B\big|\big[b,L^{-\alpha/2}\big]f_1(y)\big|^sw^s(y)\,dy\bigg)^{1/s}\bigg(\int_B w^{-s'}(y)\,dy\bigg)^{1/{s'}}\\
&\le C\|b\|_*\cdot\frac{1}{|B|}\bigg(\int_{2B}|f(y)|^pw^p(y)\,dy\bigg)^{1/p}\bigg(\int_B w^{-s'}(y)\,dy\bigg)^{1/{s'}}\\
&\le C\|b\|_*\|f\|_{L^{p,\kappa}(w^p,w^s)}\cdot\frac{w^s(2B)^{\kappa/p}}{w^s(B)^{1/s}}\\
&\le C\|b\|_*\|f\|_{L^{p,\kappa}(w^p,w^s)},
\end{split}
\end{equation*}
where the last inequality is due to our assumption $\kappa=p/s$. For the term $J_2$, since the kernel of $e^{-t_BL}$ is $p_{t_B}(y,z)$, then we may write
\begin{equation*}
\begin{split}
J_2\le&\,\frac{1}{|B|}\int_B\int_{\mathbb R^n}\big|p_{t_B}(y,z)\big|\big|\big[b,L^{-\alpha/2}\big]f_1(z)\big|dzdy\\
\le&\,\frac{1}{|B|}\int_B\int_{2B}\big|p_{t_B}(y,z)\big|\big|\big[b,L^{-\alpha/2}\big]f_1(z)\big|dzdy\\
&+\sum_{k=1}^\infty\frac{1}{|B|}\int_B\int_{2^{k+1}B\backslash2^kB}
\big|p_{t_B}(y,z)\big|\big|\big[b,L^{-\alpha/2}\big]f_1(z)\big|dzdy\\
=&\,J'_2+J''_2.
\end{split}
\end{equation*}
For any $y\in B$ and $z\in 2B$, by (1.1), we have $\big|p_{t_B}(y,z)\big|\le C\cdot(t_B)^{-n/2}$. Thus
\begin{equation*}
\begin{split}
J'_2&\le C\cdot\frac{1}{|B|}\int_B\int_{2B}\frac{1}{(t_B)^{n/2}}\big|\big[b,L^{-\alpha/2}\big]f_1(z)\big|dzdy\\
&\le C\cdot\frac{1}{|2B|}\int_{2B}\big|\big[b,L^{-\alpha/2}\big]f_1(z)\big|dz.
\end{split}
\end{equation*}
Using the same arguments as in the estimate of $J_1$, we can also deduce
\begin{equation*}
J'_2\le C\|b\|_*\|f\|_{L^{p,\kappa}(w^p,w^s)}.
\end{equation*}
As before, we note that for any $y\in B$, $z\in(2B)^c$, then $|z-y|\sim|z-x_0|$. In this case, by using (1.1) again, we get $\big|p_{t_B}(y,z)\big|\le C\cdot\frac{(t_B)^{n/2}}{|y-z|^{2n}}$. Hence
\begin{equation*}
\begin{split}
J''_2&\le C\sum_{k=1}^\infty\frac{1}{|B|}
\int_B\int_{2^{k+1}B\backslash2^kB}\frac{(t_B)^{n/2}}{|y-z|^{2n}}\big|\big[b,L^{-\alpha/2}\big]f_1(z)\big|dzdy\\
&\le C\sum_{k=1}^\infty\frac{1}{2^{kn}}\frac{1}{|2^{k+1}B|}\int_{2^{k+1}B}\big|\big[b,L^{-\alpha/2}\big]f_1(z)\big|dz.
\end{split}
\end{equation*}
Following along the same lines as before, we can also show that
\begin{equation*}
\begin{split}
\int_{2^{k+1}B}\big|\big[b,L^{-\alpha/2}\big]f_1(z)\big|dz&\le C\bigg(\int_{2B}|f(z)|^pw^p(z)\,dz\bigg)^{1/p}
\bigg(\int_{2^{k+1}B} w^{-s'}(z)\,dz\bigg)^{1/{s'}}\\
&\le C\|b\|_*\|f\|_{L^{p,\kappa}(w^p,w^s)}\frac{w^s(2B)^{\kappa/p}}{w^s(2^{k+1}B)^{1/s}}\cdot\big|2^{k+1}B\big|.
\end{split}
\end{equation*}
Consequently
\begin{equation}
J''_2\le C\|b\|_*\|f\|_{L^{p,\kappa}(w^p,w^s)}
\sum_{k=1}^\infty\frac{1}{2^{kn}}\cdot\left(\frac{w^s(2B)}{w^s(2^{k+1}B)}\right)^{1/s}.
\end{equation}
Observe that $w^s\in A_{1+s/{p'}}$, then there exists a number $r^*>1$ such that $w^s\in RH_{r^*}$. Moreover, by using Lemma 2.2 again, we get
\begin{equation}
\frac{w^s(2B)}{w^s(2^{k+1}B)}\le C\left(\frac{|2B|}{|2^{k+1}B|}\right)^{{(r^*-1)}/{r^*}}.
\end{equation}
Substituting the above inequality (4.2) into (4.1), we thus obtain
\begin{equation*}
\begin{split}
J''_2&\le C\|b\|_*\|f\|_{L^{p,\kappa}(w^p,w^s)}\sum_{k=1}^\infty\left(\frac{1}{2^{kn}}\right)^{{1+(r^*-1)}/{(sr^*)}}\\
&\le C\|b\|_*\|f\|_{L^{p,\kappa}(w^p,w^s)}.
\end{split}
\end{equation*}
Summarizing the estimates of $J'_2$ and $J''_2$ derived above, we can get
\begin{equation*}
J_2\le C\|b\|_*\|f\|_{L^{p,\kappa}(w^p,w^s)}.
\end{equation*}

In order to estimate the last term $J_3$, we need the following result given in \cite{deng} (see also \cite{duong1}).

\begin{lemma}
For $0<\alpha<n$, the difference operator $(I-e^{-tL})L^{-\alpha/2}$ has an associated kernel $\widetilde K_{\alpha,t}(y,z)$ which satisfies the following estimate
\begin{equation}
\big|\widetilde K_{\alpha,t}(y,z)\big|\le\frac{C}{|y-z|^{n-\alpha}}\frac{t}{|y-z|^2}.
\end{equation}
\end{lemma}

Hence, by the above kernel estimate (4.3) and the definition of $b\in{\dot\Lambda}_\beta(\mathbb R^n)$, we have
\begin{equation*}
\begin{split}
J_3&=\frac{1}{|B|}\int_B\left|\big(I-e^{-t_BL}\big)L^{-\alpha/2}\big([b(y)-b(\cdot)]f_2\big)(y)\right|\,dy\\
&\le\frac{1}{|B|}\int_B\int_{(2B)^c}\big|\widetilde K_{\alpha,t_B}(y,z)\big||b(y)-b(z)||f(z)|\,dz\\
&\le C\|b\|_{\dot\Lambda_\beta}\cdot\frac{1}{|B|}
\int_B\int_{(2B)^c}\frac{1}{|y-z|^{n-\alpha-\beta}}\frac{r_B^2}{|y-z|^2}|f(z)|\,dzdy\\
&\le C\|b\|_{\dot\Lambda_\beta}\sum_{k=1}^\infty\frac{1}{2^{2k}}\frac{1}{|2^{k+1}B|^{1-{(\alpha+\beta)}/n}}
\int_{2^{k+1}B}|f(z)|\,dz.
\end{split}
\end{equation*}
Since $w\in A_{p,s}$, then by using the estimate (3.4) and the fact that $\kappa=p/s$, we thus obtain
\begin{equation*}
\begin{split}
\int_{2^{k+1}B}|f(z)|\,dz&\le C\|f\|_{L^{p,\kappa}(w^p,w^s)}\big|2^{k+1}B\big|^{1/s+1/{p'}}\cdot\frac{1}{w^s(2^{k+1}B)^{1/s-\kappa/p}}\notag\\
&= C\|f\|_{L^{p,\kappa}(w^p,w^s)}\big|2^{k+1}B\big|^{1-{(\alpha+\beta)}/n}.
\end{split}
\end{equation*}
Therefore
\begin{equation*}
\begin{split}
J_3&\le C\|b\|_{\dot\Lambda_\beta}\|f\|_{L^{p,\kappa}(w^p,w^s)}\sum_{k=1}^\infty\frac{1}{2^{2k}}\\
&\le C\|b\|_{\dot\Lambda_\beta}\|f\|_{L^{p,\kappa}(w^p,w^s)}.
\end{split}
\end{equation*}
Combining the above estimates for $J_1$, $J_2$ and $J_3$ and taking the supremum over all balls $B\subseteq\mathbb R^n$, we finally conclude the proof of Theorem 1.4.
\end{proof}

\end{document}